\documentclass[10pt]{article}
\usepackage{graphicx,epsfig}
\usepackage{graphicx}
\usepackage{latexsym}
\usepackage{amssymb}
\usepackage{amsmath}
\usepackage[latin1]{inputenc}
\usepackage{layout}

\newtheorem{prop}{Proposition}
\newtheorem{lemma}{Lemma}

\newtheorem{corollary}{Corollary}
\newtheorem{theorem}{Theorem}
\newtheorem{remark}{Remark}

\def\real{{\mathord{{\rm I\kern-2.8pt R}}}}        % Fake blackboard bold R.
\def\inte{{\mathord{{\rm I\kern-2.8pt N}}}}

\def\sZZ{{\rm Z\kern-2.8ptem{}Z}}

\def\z{{\mathchoice
  {\sZZ}
  {\sZZ}
  {\rm Z\kern-0.30em{}Z}
  {\rm Z\kern-0.25em{}Z} }}
\def\sQQ{{\kern 0.27em \vrule height1.45ex width0.03em depth0em
          \kern-0.30em \rm Q}}
\def\qu{{\mathchoice
    {\sQQ}
    {\sQQ}
  {\kern 0.225em \vrule height1.05ex width0.025em depth0em \kern-0.25em \rm Q}
  {\kern 0.180em \vrule height0.78ex width0.020em depth0em \kern-0.20em \rm Q}
        }}
\def\sCC{{\kern 0.27em \vrule height1.45ex width0.03em depth0em
          \kern-0.30em \rm C}}
\def\complex{{\mathchoice
    {\sCC}
    {\sCC}
  {\kern 0.225em \vrule height1.05ex width0.025em depth0em \kern-0.25em \rm C}
  {\kern 0.180em \vrule height0.78ex width0.020em depth0em \kern-0.20em \rm C}
        }}

\newcommand{\R}{\mathbb{R}}

\newcommand{\ba}{\begin{array}}
\newcommand{\ea}{\end{array}}
\newcommand{\be}{\begin{equation}}
\newcommand{\ee}{\end{equation}}
\newcommand{\bea}{\begin{eqnarray}}
\newcommand{\eea}{\end{eqnarray}}
\newcommand{\beaa}{\begin{eqnarray*}}
\newcommand{\eeaa}{\end{eqnarray*}}

\def\z{\zeta}

\font\tenmath=msbm10 \font\sevenmath=msbm7 \font\fivemath=msbm5
\newfam\mathfam \textfont\mathfam=\tenmath
\scriptfont\mathfam=\sevenmath \scriptscriptfont\mathfam=\fivemath

\def \E{I\!\!E \,}

\def \={{\buildrel {\rm (law)} \over =}}

\def \R{{\math R}}
\def \Z{{\math Z}}

\def\L{\Lambda}

\def\L{\mathbb{L}}
\def\R{\mathbb{R}}
\def\Z{\mathbb{Z}}

\def\text{\mbox}

\def\1{{\bf 1}}
\def\q{\quad}

\def\qed{ \hfill \vrule width.25cm height.25cm depth0cm\smallskip}

\newcommand{\basa}{\begin{assumption}}
\newcommand{\easa}{\end{assumption}}

\newcommand{\bas}{\begin{assum}}
\newcommand{\eas}{\end{assum}}

\newcommand{\Cov}{\mbox{Cov}\,}
\newcommand{\cov}{\mbox{Cov}\,}
\newcommand{\Var}{\mbox{Var}\,}

\def\limiteN{\renewcommand{\arraystretch}{0.5}
\begin{array}[t]{c}\stackrel{}{\longrightarrow} \\
{\scriptstyle N\rightarrow
\infty}\end{array}\renewcommand{\arraystretch}{1}}

\def\limiteloiN{\renewcommand{\arraystretch}{0.5}
\begin{array}[t]{c}\stackrel{{\cal D}}{\longrightarrow} \\
{\scriptstyle N\rightarrow
\infty}\end{array}\renewcommand{\arraystretch}{1}}

\def\limiteasN{\renewcommand{\arraystretch}{0.5}
\begin{array}[t]{c}\stackrel{{a.s.}}{\longrightarrow} \\
{\scriptstyle N\rightarrow
\infty}\end{array}\renewcommand{\arraystretch}{1}}

\def\limiteL2{\renewcommand{\arraystretch}{0.5}
\begin{array}[t]{c}\stackrel{\L^{2}(\Omega)}{\longrightarrow} \\
{\scriptstyle N\rightarrow
\infty}\end{array}\renewcommand{\arraystretch}{1}}

\def\limiteLL{\renewcommand{\arraystretch}{0.5}
\begin{array}[t]{c}\stackrel{\L^{2}(\R^2)}{\longrightarrow} \\
{\scriptstyle N\rightarrow
\infty}\end{array}\renewcommand{\arraystretch}{1}}

\def\simloi{\renewcommand{\arraystretch}{0.5}
\begin{array}[t]{c}\stackrel{{\cal D}}{\sim} \\
{}\end{array}\renewcommand{\arraystretch}{1}}

\newcommand{\ignore}[1]{}
\begin{document}

\renewcommand{\thefootnote}{\fnsymbol{footnote}}

\title{Asymptotic behavior of the Whittle estimator for the increments of a Rosenblatt process}
\vskip1cm \author{Jean-Marc Bardet$^{1}$ \quad Ciprian Tudor\footnote{Supported by the CNCS grant PN-II-ID-PCCE-2011-2-0015. Associate member of the team Samm, Universit\'e de Panth\'eon-Sorbonne Paris 1. Partially supported by the ANR grant "Masterie" BLAN 012103.}\vspace*{0.1in}$^{2, 3}$\vspace*{0.1in}\\$^{1}$ S.A.M.M., Universit\'{e} de Paris 1 Panth\'eon-Sorbonne\\90, rue de Tolbiac, 75634, Paris, France.\vspace*{0.1in}\\$^{2}$ Laboratoire Paul Painlev\'e, Universit\'e de Lille 1 \\
 F-59655 Villeneuve d'Ascq, France.\vspace*{0.1in}\\ 
   \quad $^{3}$ 
Academy of Economical Studies\\
Piata Romana nr. 6, Sector 1, Bucharest, Romania.
}
\maketitle

\begin{abstract}
The purpose of this paper is to estimate the self-similarity index of the Rosenblatt process by using the Whittle estimator. Via chaos expansion into multiple stochastic integrals, we establish a non-central limit theorem satisfied by this estimator. We illustrate our results by numerical simulations.
\end{abstract}

\vskip0.5cm

\textbf{2000 AMS Classification Numbers: }Primary: 60G18; Secondary 60F05, 60H05, 62F12.

\vskip0.2cm


\vskip0.5cm

\section{Introduction}
The Rosenblatt process appears as limit of normalized sums of long-range dependent series (see \cite{DM}, \cite{Ta1}). In the last years, this stochastic processes has been the object of several research papers (see \cite{PiTa}, \cite{T}, \cite{TV}, \cite{VeilletteTaqqu2012b} among others).
Its analysis is motivated by the fact that the Rosenblatt process is self-similar with stationary increments and moreover its exhibits long-range dependence (or long-memory). In this sense, it shares many properties with the more known fractional Brownian motion (fBm in the sequel) except the fact that the Rosenblatt process is not a Gaussian  process. Recall that the fBm is the only Gaussian self-similar process with stationary increments.

The practical aspects of Rosenblatt process are striking: it provides  a new class of processes
from which to model long memory, self-similarity, and H\"{o}lder-regularity,
allowing significant deviation from fBm and other Gaussian processes.  The need of non-Gaussian self-similar processes in practice (for example in
hydrology) is mentioned in the paper \cite{Taqqu3} based on the study of stochastic modeling for
river-flow time series in \cite{LK}. 

The Hurst parameter $H$ characterizes all the important properties of a
Rosenblatt process, as seen above. Therefore, estimating $H$ properly is an important problem in the analysis of this process. The Hurst parameter estimation from a $N$-length path of a fBm or more generally of self-similar or long-range dependent processes, has a long history. Several statistics have been introduced and studied to this end, such as parametric estimators (maximum likelihood estimator) as well as semi-parametric estimators (spectral, variogram or wavelets based estimators). Informations on these various approaches can be found in the books of Beran \cite{B} or Doukhan {\it et al.} \cite{Doukhan}. But in the particular case of the fBm, the main results concerning this estimation were certainly obtained by \cite{fox-taq1} and \cite{dahl89}. Indeed the optimal method for estimating $H$ is obtained from the maximization of the Whittle approximation of the log-likelihood introduced by Whittle in \cite{W} (in the sequel, the Whittle estimator). This estimator shares the same asymptotic behavior than the maximum likelihood estimator (MLE), notably it is asymptotically efficient, but numerically the Whittle estimator is clearly many more 
interesting than the MLE (no need to inverse the covariance matrix). These properties also hold for long memory stationary Gaussian  processes as it was established in \cite{fox-taq1} and \cite{dahl89} under almost general conditions: the Whittle estimator is asymptotically normal with a $\sqrt N$ convergence rate.

In the case of long memory non-Gaussian time series, there exist very few results concerning the limit behavior of the Whittle estimator. In \cite{gi-sur}, the case of fourth order moment linear processes have been considered and it was proved that the Whittle estimator is still asymptotically normal with a $\sqrt N$ convergence rate. But in \cite{GT}, the cases of functionals of long memory Gaussian processes have been studied and the conclusion is different: in general, the Whittle estimator satisfies a non-central limit theorem with a non Gaussian limit distribution and a convergence rate smaller than $\sqrt N$. This is notably the case when the functional is the Hermite polynomial $H_2(x)=x^2-1$. 

Estimating the memory parameter of the Rosenblatt process appears to be a challenging problem. This is due to the fact that this process is not a Gaussian process, the explicit expression of its probability density is not know and standard techniques cannot be applied in this case. The development of new criteria, based on the Malliavin calculus and chaos expansion into multiple Wiener-It\^o integrals (see the monograph \cite{NPbook}), for the convergence of sequences of random variables recently led to new results. In \cite{TV} and \cite{CTV} the authors studied the asymptotic behavior of the  quadratic variations of the Rosenblatt process in order to obtain the asymptotic properties of an estimator for the self-similarity index. An approach based on wavelets has been also proposed in \cite{BaTu}. A common denominator of all these works is that the estimators constructed are consistent but in general their limit behavior is not Gaussian.  This not very convenient for practical aspects.

We want to put a new brick to the theory of the long-memory parameter estimation for non-Gaussian stochastic processes. We analyze the limit behavior of the Whittle estimator for the self-similarity index $H$ of the increments of a Rosenblatt process (where $0.5<H<1$). We will see that, as in the case of the estimators based on the quadratic variations (see \cite{TV}), the Whittle estimator has a non-Gaussian limit behavior and the convergence rate is $N^{1-H}$. This is due to the fact that, if one compares with the fBm case, the chaos expansion of the estimator involves a new term with a strong dependence structure that cannot be compensated by the smoothness of the estimator. This result is not totally surprising. Indeed, from \cite{Ta1}, the Rosenblatt can be obtained as the limit of $\frac 1 n \sum_{k=1}^{[nt]} H_2(X_k)$ when $n\to \infty$, where $(X_k)$ is a long memory Gaussian process and we know from \cite{GT} that the Whittle estimator of $(H_2(X_k))_k$ also satisfies the same kind of non-
central limit theorem. 

Unfortunately, our new result concerning the Whittle estimator keeps open the following question: is it possible to propose an asymptotically Gaussian estimator of the self-similarity parameter of a Rosenblatt process? However, Monte-Carlo experiments attests that the Whittle  estimator numerically provides accurate estimations, clearly better than with other well known estimators and with almost the same quality than the one obtained for the fBm when $H$ is close to $0.5$. Hence, even if it asymptotically satisfies a non-central limit theorem, the Whittle estimator is really interesting for estimating the $H$ parameter of a Rosenblatt process.\\

We organized  our paper as follows. Section 2 contains some preliminaries on multiple stochastic integrals, the Rosenblatt process and the Whittle estimator. In Section 3 we analyze the asymptotic behavior of the Whittle estimator for the self-similarity index. Finally, Section \ref{simu} contains a numerical study of the estimator and main proofs of this paper are established in Section \ref{proofs}.

\section{Preliminaries}

In this section we introduce the basic concepts  used   throughout the paper.  We will present the  the multiple stochastic integrals, the definition and the immediate properties of the Rosenblatt process  and the basic facts concerning the Whittle estimator.

To begin with we call back some elements about multiple stochastic integrals.
\subsection{Multiple stochastic integrals}
Let $W=(W_{t})_{t\in \R}$ be a classical Wiener process on a standard
probability space $\left(\Omega,{\mathcal{F}},\mathbf{P}\right)$. If
$f\in L^{2}(\R^{n})$ with $n\geq 1$ integer, we introduce the
multiple Wiener-It\^{o} integral of $f$ with respect to $W$. The
basic reference is the monograph  \cite{N}. Let $f\in
{\mathcal{S}_{n}}$ be an elementary function with $n$ variables that
can be written as $ f=\sum_{i_{1},\ldots ,i_{n}}c_{i_{1},\ldots
,i_{n}}1_{A_{i_{1}}\times \ldots \times A_{i_{n}}}$, where the
coefficients satisfy $c_{i_{1},\ldots ,i_{n}}=0$ if two indexes
$i_{k}$ and $i_{l}$ are equal and the sets $A_{i}\in
{\mathcal{B}}(\R)$ are pairwise disjoint. For  such a step function
$f$ we define
\begin{equation*}
I_{n}(f)=\sum_{i_{1},\ldots ,i_{n}}c_{i_{1},\ldots
,i_{n}}W(A_{i_{1}})\ldots W(A_{i_{n}})
\end{equation*}
where we put $W(A)=\int_{\R} 1_{A}(s)dW_{s}$. It can be seen that
the application $ I_{n}$ constructed above from ${\mathcal{S}}_{n}$
to $L^{2}(\Omega )$ is an isometry on ${\mathcal{S}}_{n}$  in the
sense
\begin{equation}
\E\left[ I_{n}(f)I_{m}(g)\right] =n!\langle
f,g\rangle_{L^{2}(T^{n})}\q\mbox{ if }m=n  \label{isom}
\end{equation}
and
\begin{equation*}
\E\left[ I_{n}(f)I_{m}(g)\right] =0\q\mbox{ if }m\not=n.
\end{equation*}
\noindent Since the set ${\mathcal{S}_{n}}$ is dense in
$L^{2}(\R^{n})$ for every $n\geq 1$ the mapping $ I_{n}$ can be
extended to an isometry from $L^{2}(\R^{n})$ to $L^{2}(\Omega)$ and
the above properties hold true for this extension.
\\
\noindent It also holds that $ I_{n}(f) = I_{n}\big( \tilde{f}\big)$,
where $\tilde{f} $ denotes the symmetrization of $f$ defined by
$$\tilde{f}(x_{1}, \ldots , x_{n}) =\frac{1}{n!}
\sum_{\sigma}f(x_{\sigma (1) }, \ldots , x_{\sigma (n) } ),$$
$\sigma$ running over all permutations of $\left\{1,...,n\right\}$.
We will need the general formula for calculating products of Wiener
chaos integrals of any orders $m,n$ for any symmetric integrands
$f\in L^{2}(\R^{m})$ and $g\in L^{2}(\R^{n})$, which is
\begin{equation}
I_{m}(f)I_{n}(g)=\sum_{\ell=0}^{m\wedge n}\ell!\binom{m}{\ell}\binom{n}{\ell}%
I_{m+n-2\ell}(f\otimes_{\ell}g) \label{product},
\end{equation}
where the contraction $f\otimes_{\ell}g$ is  defined by
\begin{eqnarray}
 &&(f\otimes_{\ell} g) ( s_{1}, \ldots, s_{m-\ell}, t_{1}, \ldots,t_{n-\ell})\nonumber\\
&& =\int_{T ^{m+n-2\ell} } f( s_{1}, \ldots, s_{m-\ell}, u_{1},
\ldots,u_{\ell})g(t_{1}, \ldots, t_{n-\ell},u_{1}, \ldots,u_{\ell})
du_{1}\ldots du_{\ell} \label{contra}
\end{eqnarray}
and by extension $f\otimes_{0} g=f\otimes g$. Note that the contraction $(f\otimes_{\ell} g) $ is an element of
$L^{2}(\R^{m+n-2\ell})$ but it is not necessarily symmetric. We will
denote its symmetrization by $(f \tilde{\otimes}_{\ell} g)$.
\\

\subsection{The Rosenblatt process}
Recall that a fBm $B=(B_t)_{t\geq 0}$ with Hurst parameter $H\in (\frac 1 2 ,1)$ and parameter $C>0$ is a centered Gaussian process with covariance function
\begin{equation}\label{covFBM}
\Cov(B_t,B_s)=\E(B_tB_s)=\frac {C^{2}} 2\big(s^{2H}+t^{2H}-|t-s|^{2H}\big),\quad
s,t\in[0,\infty),
\end{equation}
with $C=\Var (B_1)$. It is the only normalized Gaussian $H$-self-similar process with stationary increments.
The fBm admits the following moving average representation: for every $t\geq 0$ and for every $H>\frac{1}{2}$
\begin{equation}
\label{repB}
B _{t} =Cc_{B}(H)\int_{\mathbb{R}} \left( \int_{0} ^{t} (u-y) _{+} ^{H-\frac{3}{2}} du \right) dW(y).
\end{equation}
where $(W(y), y\in \mathbb{R} )$ is a standard Brownian motion with time interval $\mathbb{R}$ and $c_{B}(H)$  is a strictly positive constant that ensures
that $E( B _{t}) ^{2}= C^{2}t^{2H}$  for every $t\geq 0$.

The Rosenblatt process is related to the fractional Brownian motion. It shares many properties of the fBm. For instance, it has the same covariance function (\ref{covFBM}) as the fBm, it is $H$-self-similar and it has stationary increments. It has the same order of the H\"older regularity of its sample paths as the fBm (that is, the Rosenblatt process, as well as the fBm, are H\"older continuous of order $\delta $ with $0<\delta <H$). There are also some differences with respect to the fBm. One of them, is that it is defined only for the self-similarity index $H>\frac{1}{2}$ and another difference, more important, is that it is not Gaussian. It can be expressed as a double multiple integral with respect to the Wiener process and therefore it is an element of the second Wiener chaos. More exactly, a Rosenblatt process $Z=(Z_t)_{t\geq 0}$ with self-similarity order $H\in (\frac{1}{2}, 1)$ and parameter $C>0$
is defined by
\begin{equation}
\label{Z}
Z_{t}= C c_{Z}(H)\int_{\mathbb{R}}\int_{\mathbb{R}}\left( \int_{0} ^{t} (u-y_{1}) _{+} ^{\frac{H}{2}-1}(u-y_{2}) _{+} ^{\frac{H}{2}-1}du \right) dW(y_{1}) dW(y_{2})
\end{equation}
where $ c_{Z}(H)$ is a strictly positive constant such that $E(Z_{1})^{2}= C^{2}$, {\it i.e.}
\begin{eqnarray}\label{CZ}
 c^2_Z(H)=\frac {2H(2H-1)}{\beta^2(1-H,\frac H 2)} \qquad \mbox{see e.g. \cite{T}}.
\end{eqnarray}
We will call a (standard) Rosenblatt random variable every random variable that has the same law as $Z_{1}$ with parameter $C=1$. In the sequel we will use the kernel $L_H$ of the Rosenblatt process defined by
\begin{equation}
\label{LH}
L^{H}(y_1,y_2)=  \int_{0} ^{t} (u-y_{1}) _{+} ^{\frac{H}{2}-1}(u-y_{2}) _{+} ^{\frac{H}{2}-1}du \qquad\mbox{for $(y_1,y_2)\in \R^2$.}
\end{equation}

As the second order properties of a Rosenblatt process are the same than the ones of the corresponding fBm, the process $Y=(Y_t)_{t\geq 0}$ of the increments  of a Rosenblatt process, defined by
$$
Y_t=Z_{t+1}-Z_{t}\qquad \mbox{for $t\geq 0$}
$$
with $Z_0=0$ by definition, is a long memory stationary process with covariogram
\begin{equation}\label{covY}
r_{H,C}(t):=\cov( Y_0,\, Y_t)=\frac {C^2} 2 \, \big ( |t+1|^{2H}+|t-1|^{2H}-2 |t|^{2H} \big )\quad\mbox{for $t\in \R$},
\end{equation}
and a spectral density $f_{H,C}$ defined for $\lambda\in [-\pi,\pi]$ by:
\begin{equation}\label{covfY}
f_{H,C}(\lambda):=\frac 1 {2\pi} \, \sum_{k \in \Z} r(k)\, e^{ik\lambda}= \frac {C^2 \, H\Gamma(2H)\sin(\pi H)}
{2\pi } \,  (1-\cos \lambda) \sum_{k\in \Z} |\lambda+2k\pi|^{-1-2H},
\end{equation}
since $\int_\R (1-\cos x)|x|^{-1-2H}dx =2\pi \big(H\Gamma(2H)\sin(\pi H)\big)^{-1}$ (see Sinai, 1976, or Fox and Taqqu, 1986).

\subsection{The Whittle estimator}
Our purpose is to study the asymptotic properties of the Whittle estimator of parameters $H$ and $C$  computed from a sample $(Z_1,\cdots,Z_N)$ of a Rosenblatt process $Z=(Z_t)_{t\geq 0}$ .  Let us briefly introduce the Whittle estimator. The first step is to 
define the periodogram of the process $Y$ of the increments of $Z$:
\begin{equation}\label{perio}
\widehat I_N(\lambda) =\frac1{2\pi \, N} \Big|\sum_{k=0}^{N-1}
Y_k\,e^{-ik\lambda}\Big|^2%=\frac1{2\pi \, N} \Big|e^{i\lambda N}Z_N+ \sum_{k=1}^{N-1} Z_k(1-e^{i\lambda})\,e^{-ik\lambda}\Big|^2,~~~\mbox{for}~\lambda \in [-\pi,\pi].
\end{equation}
Now, let $g:\R \rightarrow \R$ be a $2\pi$-periodic function such
that $g\in\L^2([-\pi,\pi])$ and define
\begin{eqnarray*}
\widehat J_N(g)&=& \int_{-\pi}^\pi g(\lambda)\widehat I_n(\lambda)\,
d\lambda,~~~\mbox{the integrated periodogram of}~Y
\\
\mbox{and}~~J(g)&=& \int_{-\pi}^\pi g(\lambda)f_{H,C}(\lambda)\, d\lambda,
\end{eqnarray*}
with $f_{H,C}$ denotes the spectral density of $Y$ defined in  (\ref{covfY}). From (\ref{covY}) and (\ref{perio}), we also have
\begin{equation}\label{INbis}
\displaystyle{\widehat I_N(\lambda) = \frac 1 {2 \pi}
\sum_{|k|<N}  \widehat r_N(k)e^{-ik\lambda}~~ \mbox{with}~~\widehat
r_N(k)=\frac 1 N \sum_{j=1}^{N-|k|}Y_jY_{j+|k|}},
\end{equation}
which is a biased estimate of $r_{H,C}(k)$ (see the next section). Thus, the periodogram
$\widehat I_N(\lambda)$ could be a natural estimator of the spectral density;
unfortunately it is not a consistent estimator. However, once
integrated with respect to some $\L^2$ function, its behavior
becomes quite smoother and can allow an estimation of the spectral
density. The Whittle's contrast is a special case of the integrated periodogram $\widehat J_N(g)$. Indeed, this contrast can be written (see \cite{W})
$$
\widehat U_N(H,C)=\int_{-\pi}^\pi \Big ( \log(f_{H,C}(\lambda))+ \widehat J_N(1/f_{H,C}(\lambda))\Big )d\lambda.
$$
The Whittle estimator is thus:
$$(\widehat H_N,\widehat C_N):=\mbox{Argmin}_{(H,C)\in ( 1/2,1)\times (0,\infty)}   \widehat U_N(H,C).
$$
But using a classical renormalization (see Fox and Taqqu, 1986 for instance),
the spectral density $f_{H,C}$ of $Y$ can be decomposed as:
\begin{eqnarray}\label{normalized}
f_{H,C}(\lambda)&=&\sigma^2\, g_H(\lambda)\quad \mbox{with} \left \{ \begin{array}{lll}
  \sigma^2&=& \displaystyle \frac 1 {2\pi} \, \frac {C^2} {a_H}  H\Gamma(2H)\sin(\pi H)     \vspace{0.2cm}\\
  g_H(\lambda)&=&  \displaystyle a_H (1-\cos \lambda) \sum_{k\in \Z} |\lambda+2k\pi|^{-1-2H}  \\
  \end{array} \right .
\end{eqnarray}
for all $\lambda\in [-\pi,0)\cup (0,\pi]$ with $\displaystyle a_H:=\exp \Big [ -\int_{-\pi}^\pi \log\Big ( (1-\cos t) \sum_{k\in \Z} \big | t +2k\pi \big |^{-1-2H}\Big )dt\Big ]$.\\
Then
\begin{equation}
\label{io}
\int_{-\pi}^\pi \log ( g_H(\lambda))d\lambda=0
\end{equation}
 for $H\in (1/2,1)$ and the minimization of $\widehat U_N(H,C)$ can be write again as a minimization in $(H,\sigma^2)$ and this implies that
\begin{eqnarray}\label{defestim}
\widehat{H}_N &=& \mbox{Argmin} _{H \in (1/2,1)} \left \{\widehat
J_N(1/g_H) \right \}=\mbox{Argmin} _{H \in (1/2,1)}\Big \{  \int_{-\pi}^\pi \frac
{\widehat I_N(\lambda)}{g_H(\lambda)}\, d\lambda \Big \}\\
\label{defestim2}~~\mbox{and}~~~\widehat C_N&=&\Big (2\pi \, \frac {a_{\widehat H_N} \, \widehat{\sigma}^2_N }{{\widehat H_N}\Gamma(2{\widehat H_N})\sin(\pi {\widehat H_N})  }\Big )^{1/2} \qquad \mbox{with}\qquad  \widehat{\sigma}^2_N :=\frac 1
{2\pi}\widehat J_N (1/g_{\widehat H_N}).
\end{eqnarray}
\begin{remark}
However, for practical use, these definitions $\widehat{H}_N$
and $\widehat{\sigma}^2_N$ have to be modified. Indeed, the
assumption that the process has zero mean is unrealistic.
Moreover, the integrals defining the estimators has to be replaced
by their approximations by a Riemann sum. Thus, define:
\begin{eqnarray*}\label{estim2}
\widetilde{H}_N~&=& \mbox{Argmin} _{H \in (1/2,1)}\Big \{
\frac {2\pi} N \, \sum_{k=1}^N  \frac {\widetilde I_N(\pi
k/N)} {g_H(\pi k/N)} \Big
\}~~~\mbox{and}~~~\widetilde{\sigma}^2_N =\frac 1  N \, \sum_{k=1}^N \frac {\tilde I_N(\pi k/N)}{g_{\widetilde{H}_N}(\pi
k/N)},\\
&& \mbox{with}~~\widetilde I_N(\lambda)=\frac1{2\pi \, N}
\Big|\sum_{k=0}^{N-1}
(Y_k-\overline{Y}_N)e^{-ik\lambda}\Big|^2~~\mbox{for}~\lambda \in
[-\pi,\pi]~~\mbox{and}~~\overline{Y}_N=\frac 1 N \, \sum_{k=0}^{N-1}
Y_k.
\end{eqnarray*}
\end{remark}

\section{Limit theorems for the Whittle estimator of parameters of a Rosenblatt process}

In this section, we prove limit theorems for the integrated periodogram $\widehat J_{N}(g)$ for $g \in \L^2[-\pi,\pi])$ satisfying certain conditions. Applied to  $g=1/g_H$, we will prove the almost-sure convergence of $\widehat H_N$. Applied to $g=\frac {\partial}{\partial H} (1/g_H)$  and then with a classical Taylor expansion also using the case $g=\frac {\partial^2}{\partial^2 H} (1/g_H)$ this will provide a non-central limit theorem satisfied by $\widehat H_N$. \\
~\\
Therefore, the main point is to obtain limit theorems for $\widehat J_{N}(g)$ . For this, we use a (now) standard approach based on the chaotic decomposition of the random variable $\widehat J_{N}$ into a sum of multiple stochastic integrals. Since the Rosenblatt process at fixed time is a multiple integral of order $2$, the product $Y_{j} Y_{j+k}$ is a product of two double multiple integrals which can be expressed as a the sum of a multiple integral of order $4$, a multiple integral of order $2$ and a deterministic function, which is the expectation of $Y_{j} Y_{j+k}$. This decomposition is transferred to the integrated periodogram $\widehat J_{N}$, which will be written as sum of two multiple integrals (one of order 2, one of order 4) plus its expectation. What we  do next, is to analyze  these three terms that compose $\widehat J_{N}$. We will see that, as in the case of the variation statistic of the Rosenblatt process (see \cite{TV}) or of wavelet statistic (see \cite{BaTu}), the dominant term is the
one in the second Wiener chaos,
which will give the asymptotic behavior of $\widehat J_{N}$. A detailed study of this term shows that it converges to a Rosenblatt random variable. 
\subsection{Chaos decomposition of the integrated periodogram }
The purpose of this part is to provide the asymptotic behavior of a sequence $\widehat{J} _{N} (g)$ using its chaos decomposition. We start first with the analysis of $\widehat{r} _{N} (k)$ with fixed $k \in \mathbb{Z}$.
%Let us assume for the moment that $k\geq 0$ because the case $k<0$ can be treated in a symmetric way.
%In this case
%\begin{equation*}
%\widehat{r} _{N} (k)= \frac{1}{N}   \sum_{j=0}^{N-|k|-1} Y_{j} Y_{j+k}.
%\end{equation*}
We will study the convergence of $\widehat{r} _{N} (k)$ to the covariance function $r_{H,C}(k)$. % for fixed $k\geq 0$.
We first observe that
\begin{equation*}
\E \widehat{r} _{N} (k) =\frac{1}{N}  \sum_{j=0}^{N-|k|-1} \E Y_{j} Y_{j+|k|} = \frac{N-|k|} {N} r_{H,C}(k)
\end{equation*}
and this converges to $r_{H,C}(k)$ when $N\to \infty$: $\widehat{r} _{N} (k)$ is an asymptotically unbiased estimator for $r_{H,C}(k)$. On the other hand, for every $j\geq 1$  we can write the increment of the Rosenblatt process as
\begin{equation*}
Y_{j}= Z_{j+1}-Z_{j} = I_{2} (\Delta L_{j} ^{H} )
\end{equation*}
where $L^{H}$ is the kernel of the Rosenblatt process  (\ref{LH}) and we denoted by $\Delta L^{H} _{j} $ the two-variables kernel
\begin{eqnarray*}
\Delta L ^{H} _{j} (y_{1}, y_{2})&=&C\, c_{Z}(H) \int_{j} ^{j+1}(u-y_{1}) _{+} ^{\frac{H}{2}-1} (u-y_{2}) _{+} ^{\frac{H}{2}-1} du.
\end{eqnarray*}
In the sequel we will simply denote $L^{H}:=L$. Then, for every $|k|<N$, by the product formula (\ref{product})
\begin{eqnarray*}
\widehat r_{N}(k)-r_{H,C}(k)&=& \frac{1}{N} \sum_{j=0}^{N-|k|-1} I_{2} (\Delta L_{j}) I_{2}(\Delta L_{j+|k|} ) \\
&=& \frac{1}{N} \sum_{j=0}^{N-|k|-1} \big (I_{4} (\Delta L _{j}\otimes \Delta L_{j+|k|}) + 4 \, I_{2} (\Delta L _{j} \otimes _{1} \Delta L _{j+|k|})\big )\\
&&+ \E \widehat{r} _{N}(k) -r_{H,C}(k).
\end{eqnarray*}
Therefore, for $g \in\L^2([-\pi,\pi])$,
\begin{eqnarray}
\label{JN} \widehat J_{N}(g)- \E \widehat J_{N}(g)&\hspace{-0.3cm}=&\hspace{-0.3cm}\frac{1}{2\pi \, N}\int_{-\pi} ^{\pi }\hspace{-0.2cm} d\lambda g(\lambda ) \sum_{|k|<N} \sum_{j=0}^{N-|k|-1}  e^{ik\lambda } \ \big (I_{4} (\Delta L _{j} \otimes \Delta _{j+|k|}L) +4 \,  I_{2} (\Delta L _{j} \otimes _{1} \Delta L _{j+|k|})\big )~~~\\
\nonumber &\hspace{-0.3cm}=& \hspace{-0.3cm}\widehat T_{2,N}(g) + \widehat T_{2,N}(g),\\
\label{T4} \mbox{where}\qquad \widehat T_{4, N}(g) &\hspace{-0.3cm}:=&\hspace{-0.3cm}\frac{1}{2\pi \, N}\int_{-\pi} ^{\pi }\hspace{-0.2cm} d\lambda g(\lambda ) \sum_{|k|<N}  \sum_{j=0}^{N-|k|-1}  e^{ik\lambda } \ I_4 \Big (\Delta L _{j} (y_{1},y_{2}) \Delta L _{j+|k|} (y_{1},y_{2}) \Big ) \\
\label{T2} \mbox{and}\qquad  \widehat T_{2, N}(g) &\hspace{-0.3cm}:=&\hspace{-0.3cm}\frac{2}{\pi \, N}\int_{-\pi} ^{\pi }\hspace{-0.2cm} d\lambda g(\lambda ) \sum_{|k|<N} \sum_{j=0}^{N-|k|-1}  e^{ik\lambda } \ I_2 \Big ( \int_{\mathbb{R} } dx \Delta L _{j} (y_{1}, x) \Delta L _{j+|k|} (y_{2}, x) \Big ).
\end{eqnarray}
The notation $\widehat T_{2,N}(g), ~\widehat T_{4,N}(g) $  suggests that the sequence $\widehat T_{4,N}(g) $ belongs to the fourth Wiener chaos while $\widehat T_{2,N}(g) $ belongs to the second Wiener chaos.

\subsection{Asymptotic behavior of the integrated periodogram}
In order to study the asymptotic behavior of $\widehat J_{N}(g)$ let us first specify our assumptions on the  function $g$. For $H\in (1/2,1)$ and for $g:[-\pi,\pi] \mapsto \R$, we introduce\\
~\\
{\bf Assumption A$(H)$}: $g \in \L^2([-\pi,\pi])$ and  for any $\delta>0$, there exists $C$ such that 
$$
|g(\lambda)| \leq C \, |\lambda|^{2H-1-\delta}\qquad \mbox{for all $\lambda\in [-\pi,0)\cup (0,\pi]$.}
$$
A very useful point for us is the following (already proved in \cite{FT}, p. 531):
\begin{lemma}
Let $H\in (1/2,1)$ and let  $g_H$ be given by   (\ref{covfY}). Then, $1/g_H$, $\frac{\partial}{ \partial H} (1/g_H)$ and $\frac{\partial^2}{ \partial H^2} (1/g_H)$ satisfy Assumption A$(H)$.
\end{lemma}

The limit in distribution of the sequence $\widehat J_{N}(g)- \E \widehat J_{N}(g)$ as $N\to \infty$ depends on the asymptotic behavior of the two terms $\widehat T_{2,N}(g)$ and  $\widehat T_{4,N}(g)$ above. First, let us analyze the term in the second Wiener chaos, denoted by $\widehat T_{2,N}(g)$, of the decomposition $\widehat J_{N}(g)$  when the function $g$ satisfies Assumption A$(H)$. Its asymptotic  behavior is described by the following  result.
\begin{prop}\label{p2}
For every $H\in (\frac{1}{2}, 1)$ and all function  $g$ satisfying Assumption A$(H)$,
\begin{equation*}
N^{1-H} \, \widehat T_{2,N}(g) \limiteloiN  8  \, \sqrt{ \frac{2(2H-1)}{H(1-H)^2}} \,\Big ( \int_{-\pi} ^{\pi }  g(\lambda )f_{(H+1)/2,C}(\lambda) \, d\lambda \Big )\,  R,
\end{equation*}
with $R$ a standard Rosenblatt random variable $(\E R^2=1$).
\end{prop}

Under the normalization of $\widehat T_{2,N}(g)$ the summand $\widehat T_{4,N}(g)$ goes to zero:
\begin{prop}\label{p3}
For every $H\in (\frac{1}{2}, 1)$ and all function $g$ satisfying Assumption A$(H)$,
\begin{equation*}
N^{1-H} \, \widehat T_{4,N}(g) \limiteL2 0.
\end{equation*}

\end{prop}

A non-central limit theorem satisfied by $\widehat J_N(g)$ is thus the consequence of both the previous Propositions \ref{p2} and \ref{p3}:

\begin{prop}\label{p7}
For every $H\in (\frac{1}{2}, 1)$ and all function $g$ satisfying Assumption A$(H)$,
\begin{equation*}
N^{1-H} \,\big ( \widehat J_N(g) -\E \widehat J_N(g) \big )  \limiteloiN  8  \, \sqrt{ \frac{2(2H-1)}{H(1-H)^2}} \,\Big ( \int_{-\pi} ^{\pi }  g(\lambda )f_{(H+1)/2,C}(\lambda) \, d\lambda \Big )\,  R,
\end{equation*}
with $R$ a standard Rosenblatt random variable $(\E R^2=1$).
\end{prop}

The above Proposition \ref{p7} holds for any function $g$ that satisfies Assumption A(H). The next result shows that, when $g= \frac{1}{\partial H}\left( \frac{1}{g_{H}}\right)$ (with $g_{H} $ defined in (\ref{normalized})) under the normalization of $\widehat T_{2,N}(g)$ the deterministic term in the chaos expansion of $\widehat J_{N} (g)$ converges to zero.

\begin{prop}\label{A4}
With $g_{H}$ given by (\ref{covfY}), we have
\begin{equation*}
N^{1-H} \, \E \widehat J_{N} \Big (\frac{\partial}{ \partial H} \big ( \frac 1 {g_{H}} \big ) \Big ) \limiteN 0.
\end{equation*}
\end{prop}

\vskip0.2cm

It is also possible to show the  following almost sure limit theorem for the sequence $\widehat J_{N}(g)$ when $g$ satisfies Assumption A(H).

\begin{prop}\label{p4}
For every $H\in (\frac{1}{2}, 1)$ and all function  $g$ satisfying Assumption A$(H)$, we have
\begin{equation*}
 \widehat J_{N}(g) \limiteasN \sigma ^{2}\int_{-\pi } ^{\pi} g( \lambda) \, g_{H} (\lambda) \, d\lambda.
\end{equation*}
\end{prop}

Now, we can state our main result.
\begin{theorem} \label{asymptoH}
Let $\widehat{H} _{N}$ be  defined by  (\ref{defestim}). Then
\begin{equation*}
\widehat H_N \limiteasN H \qquad \mbox{and}\qquad 
N^{1-H} \, \big( \widehat{H}_{N} -H\big ) \limiteloiN   \gamma(H)\,  R
\end{equation*}
where $R$ is a standard Rosenblatt random variable (with $\E R^2=1$) and $\gamma(H)$ is defined by:
\begin{eqnarray}\label{aH}
 \gamma(H):=  16  \pi \, \sqrt{ \frac{2(2H-1)}{H(1+H)^2}} \,\Big ( \int_{-\pi} ^{\pi } \frac {f_{(H+1)/2,1}(\lambda)}{g_H(\lambda)} \, d\lambda \Big )\,\Big (\int_{-\pi} ^{\pi} f_{H,1} (\lambda) \frac{ \partial ^{2} }{\partial H^{2} } \left( \frac{1}{g_{H} (\lambda)} \right) d\lambda \Big ) ^{-1}
\end{eqnarray}

\end{theorem}
It is also possible to provide the strong convergence of $\widehat C_N$ (defined in \ref{defestim2}) to $C$:
\begin{corollary}\label{sigma}
Under the assumptions of Theorem \ref{asymptoH}, 
$$
\widehat C_N \limiteasN C  \qquad \mbox{and}\qquad 
N^{1-H} \, \big( \widehat{C}_{N} -C\big ) \limiteloiN   C \, \rho(H)\,  R,
$$
where $R$ is a standard Rosenblatt random variable (with $\E R^2=1$)  and $\rho(H)$ is defined in \eqref{eq6}.
\end{corollary}

\section{Monte-Carlo experiments} \label{simu}
\begin{center}
\begin{figure}
\[
\includegraphics[width=12 cm,height=5 cm]{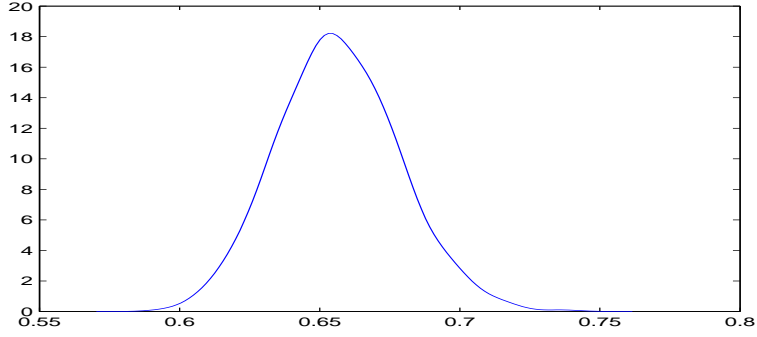}
\]
\caption{\it FFT estimation (Silverman's method) of the density of the limit of $\widehat H_N$ for $H=0.65$, $N=5000$  from $1000$ independent replications of Rosenblatt process paths.}\label{Fig1}
\end{figure}
\end{center}

We generated $200$ paths of Rosenblatt processes for several values of $H$ ($=0.55$, $0.65,\cdots,0.95$) and $N$ $(=1000, \, 5000$ and $20000$). These paths are obtained from the algorithm already used in \cite{BaTu} and deduced from the asymptotic behavior of $n^{-H} \sum_{i=1}^{[nt]} (Y_i^2-1)$ when $n\to \infty$ where $(Y_i)$ is a sequence of centered and normalized LRD processes with memory parameter $d=H/2$ (typically FARIMA$(0,H/2,0)$ processes). Note that the generator of Rosenblatt process paths as well as the computations of the estimators used in this section are available on the website {\tt http://samm.univ-paris1.fr/-Jean-Marc-Bardet} with a free access on (in Matlab language).  \\
We computed the Whittle estimator $\widehat H_N$ of $H$.  An example of the estimation of the probability density function of $\widehat H_N$ provided by the Silverman's nonparametric procedure is given in Figure \ref{Fig1}. This estimated probability density function appears could first appear as Gaussian density function but it is slightly asymmetric as the Rosenblatt density function should be when $H=0.65$ (see for instance Figure 1 in \cite{VeilletteTaqqu2012b}). \\
The Whittle estimator $\widehat H_N$ of $H$ applied to paths of Rosenblatt processes is also compared to both other estimators:
\begin{itemize}
 \item The extended local Whittle estimator $\widehat H_{ADG}$ defined in \cite{ADG} from the seminal paper \cite{Ro}. The trimming parameter is chosen as $m=N^{0.65}$ (this is not an adaptive estimator) following the numerical recommendations of   \cite{ADG}.
\item The wavelet estimator $\widehat H_{Wa}$ as it was already defined in \cite{BaTu}. 
\end{itemize}
Note that both these estimators are semi-parametric estimators (while $\widehat H_N$ is a parametric estimator) and thus can be applied to other processes than Rosenblatt processes. However the asymptotic study of $\widehat H_{ADG}$ for Rosenblatt was not already done. The results are detailed in Table \ref{Table1}. \\

\begin{table}[t]
{\footnotesize
\begin{center}
\begin{tabular}{|c|c|c|c||c|c|}
\hline\hline
  $N=1000$  & $H=0.55$ & $H=0.65$ & $H=0.75$  &$H=0.85$ & $H=0.95$  \\
%&              &$\rho=0$  &$\rho=0$  &$\rho=0$  &$\rho=0$  &$\rho=1$  &$\rho=1$  &$\rho=1$  &$\rho=1$ \\
\hline \hline
mean $\widehat H_N$    &0.570 &   0.653 &    0.736   & 0.815 &    0.917 
\\
std $\widehat H_{N}$    &0.030 &   0.041 &    0.047 &    0.053 &    0.050  
\\
\hline
mean $\widehat H_{ADG}$    &0.570 &    0.634 &    0.708 &    0.795 &    0.906 
\\
std $\widehat H_{ADG}$    &0.072 &    0.084 &    0.094&    0.105  &  0.102
\\
\hline
mean $\widehat H_{Wa}$    & 0.499 &    0.542 &   0.619 &    0.685 &    0.766
\\
std $\widehat H_{Wa}$    & 0.104 &   0.116 &    0.115 &    0.129 &   0.119
\\
\hline

\end{tabular}\\
~\\\vspace{3mm}
\begin{tabular}{|c|c|c|c||c|c|}
\hline\hline
  $N=5000$  & $H=0.55$ & $H=0.65$ & $H=0.75$  &$H=0.85$ & $H=0.95$  \\
%&              &$\rho=0$  &$\rho=0$  &$\rho=0$  &$\rho=0$  &$\rho=1$  &$\rho=1$  &$\rho=1$  &$\rho=1$ \\
\hline \hline
mean $\widehat H_N$    & 0.582 &    0.655 &    0.743 &    0.837 &    0.929
\\
std $\widehat H_{N}$    & 0.014 &    0.019&    0.029 &    0.033 &    0.035
\\
\hline
mean $\widehat H_{ADG}$    &0.575 &   0.627 &    0.723 &    0.824 &    0.919 
\\
std $\widehat H_{ADG}$    &0.041 &    0.052 &    0.062  &  0.067 &    0.072 
\\
\hline
mean $\widehat H_{Wa}$    &0.550 &  0.610 &    0.698 &    0.800 &    0.891 
\\
std $\widehat H_{Wa}$    &0.055 &    0.062 &    0.072 &    0.079 &    0.075 
\\
\hline
\end{tabular}\\
~\\\vspace{3mm}
\begin{tabular}{|c|c|c|c||c|c|}
\hline\hline
  $N=20000$  & $H=0.55$ & $H=0.65$ & $H=0.75$  &$H=0.85$ & $H=0.95$  \\
%&              &$\rho=0$  &$\rho=0$  &$\rho=0$  &$\rho=0$  &$\rho=1$  &$\rho=1$  &$\rho=1$  &$\rho=1$ \\
\hline \hline
mean $\widehat H_N$    & 0.571   & 0.656 &    0.746  &0.847 	&    0.937
\\
std $\widehat H_{N}$    &0.008   &0.015  &  0.020   & 0.025 &    0.025 
\\
\hline
mean $\widehat H_{ADG}$    &0.563 & 0.637 &   0.734   & 0.838 &  0.931
\\
std $\widehat H_{ADG}$    &0.033 & 0.040 &   0.052   & 0.040 &   0.045
\\
\hline
mean $\widehat H_{Wa}$    &0.569  & 0.630 &    0.728  &0.838 &    0.931
\\
std $\widehat H_{Wa}$    &0.040   & 0.039 &    0.052   & 0.053 &   0.042
\\
\hline
\end{tabular}
\end{center}
}
\caption{{\small\label{Table1} : Comparison between the Whittle estimator $\widehat H_N$  and other famous semiparametric estimators of $H$ (extended local Whittle estimator $\widehat H_{ADG}$ and wavelet estimator $\widehat H_{Wa}$) applied to $100$ paths of Rosenblatt processes with several $H$ and $N$ values}
}
\end{table}
{\bf Conclusions of the simulations:}  The results obtained with the Whittle estimator $\widehat H_N$ are very convincing and fit the limit theorem we established. Indeed the estimator $\widehat H_N$ seems to be asymptotically unbiased and its standard deviation is depending on $H$ (the larger $H$ the larger the standard deviation). \\
For giving a comparison, it could be interesting to compare these results with those obtained for fractional Brownian motions (fBms) instead of Rosenblatt processes (the computation of the Whittle estimator is exactly the same for both those processes). In additional simulations, we obtained the following results: the standard deviation $\overline s_N$ of  $\widehat H_N$ applied to fBms almost not depends on $H$. Hence for almost all $H \in ]0.5,1[$ and $N=1000$,   $\overline s_N\simeq 0.02$, for $N=5000$, $\overline s_N \simeq 0.01$ and for $N=20000$, $\overline s_N\simeq 0.005$. As a consequence, from Table \ref{Table1}, 
\begin{itemize}
 \item when $H=0.55$, the standard deviation of $\widehat H_N$ for Rosenblatt process is close to the one obtained for fBm. In terms of theoretical convergence rate this is not surprising since the convergence rate of $\widehat H_N$ is  $N^{0.45}$ in case of Rosenblatt process while it is $N^{0.5}$ in case of fBm. 
\item when $H=0.95$, the standard deviation of $\widehat H_N$ for Rosenblatt process is dramatically larger than the one obtained for fBm. Once again, this is not surprising since the convergence rate of $\widehat H_N$ is  $N^{0.05}$ in case of Rosenblatt process while it is still $N^{0.5}$ in case of fBm. 
\end{itemize}
Also from Table \ref{Table1}, we can compare the convergence rates of $\widehat H_N$ and both the semiparametric estimators $\widehat H_{ADG}$ and $\widehat H_{Wa}$. The standard deviation of $\widehat H_N$ is almost the half of those of $\widehat H_{ADG}$ and $\widehat H_{Wa}$. The results are really convincing and show the good accuracy of the Whittle estimator for estimating $H$. However, to be fair, we have to underline that $\widehat H_N$ is a parametric estimator typically appropriated to Rosenblatt (or fBm) increment processes while $\widehat H_{ADG}$ and $\widehat H_{Wa}$ are semiparametric estimators which can be applied to general classes of long-memory processes. 
\begin{remark}
As we improved the generator of Rosenblatt processes, the results obtained with $\widehat H_{Wa}$ are a little better than those obtained in \cite{BaTu}, especially  when $H$ is close to $1$.
\end{remark}

\section{Proofs}\label{proofs}
The following technical lemma will be needed several times in the sequel.

\begin{lemma}\label{l1}
Let $g$ satisfy Assumption A$(H)$ with $H\in (1/2,1)$. Then, for any $\delta>0$ and $k\in \Z$,
$$
\int_{-\pi } ^{\pi} g(\lambda)\, e^{ik\lambda }  d\lambda =O\big ( (1+|k|)^{-2H+\delta}\big ).
$$
%where the symbol $\sim$ means that the two sides have the same limit as $k\to \infty$.
\end{lemma}
{\bf Proof: } This have been stated and proved in \cite{FT}, Lemma 5. \qed \\
~\\
{\bf Proof of Proposition \ref{p2}: } First, using the definition of the contraction $\otimes _{1}$ (see (\ref{contra})),
\begin{eqnarray*}
 (\Delta L _{j} \otimes _{1} \Delta L _{j+|k|})(y_{1}, y_{2})
&=& \int_{\mathbb{R} } dx \Delta L _{j} (y_{1}, x) \Delta L _{j+|k|} (y_{2}, x) \\
&=& C^2 \, c^2_{Z}(H)  \int_{\mathbb{R}} dx \int_{j} ^{j+1} du(u-y_{1}) _{+} ^{\frac{H}{2}-1} (u-x) _{+} ^{\frac{H}{2}-1}
\int_{j+|k|} ^{j+|k|+1} dv(v-y_{2}) _{+} ^{\frac{H}{2}-1} (v-x) _{+} ^{\frac{H}{2}-1}\\
&=& C^2 \,  c^2_{Z}(H)\beta (1-H, \frac{H}{2} )\int_{j} ^{j+1} du\int_{j+|k|} ^{j+|k|+1} dv(u-y_{1}) _{+} ^{\frac{H}{2}-1} (v-y_{2}) _{+} ^{\frac{H}{2}-1}\vert u-v\vert ^{H-1}
\end{eqnarray*}
where we changed the order of integration and we used the identity:  for $(a,b) \in (-1,0)^2$ and $a+b<-1$, for $(u,v)\in \R^2$ with $u\neq v$,
\begin{equation}
\label{i}
\int_{-\infty} ^{u\wedge v} (u-y) ^{a} (v-y) ^{b} dy = \beta ( -1-a-b, a+1) \, |u-v|^{a+b+1}.
\end{equation}
Then
$$
\widehat T_{2, N}(g)=C^2 \,  \frac{4H(2H-1) }{\beta (1-H, \frac{H}{2} )\pi \, N}\int_{-\pi} ^{\pi }\hspace{-0.3cm} d\lambda \,  g(\lambda )\sum_{|k|<N}  \sum_{j=0}^{N-|k|-1} e^{ik\lambda } \ I_2 \Big (\int_{j} ^{j+1} du\int_{j+|k|} ^{j+|k|+1} dv(u-y_{1}) _{+} ^{\frac{H}{2}-1} (v-y_{2}) _{+} ^{\frac{H}{2}-1}\vert u-v\vert ^{H-1} \Big ).
$$
By Lemma \ref{A22}, the sequence $(N^{1-H} \widehat T_{2, N}(g))_N$   has the same limit in $\L^{2}(\Omega)$ as the sequence $(N^{1-H}\widehat T'_{2, N}(g) )_N$ where
\begin{multline*}
\widehat T'_{2, N}(g)= C^2 \,  \frac{4H(2H-1) }{\beta (1-H, \frac{H}{2} )} \, \frac 1 {\pi \, N}\int_{-\pi} ^{\pi }\hspace{-0.3cm} d\lambda \,  g(\lambda )\sum_{|k|<N} \sum_{j=0}^{N-|k|-1} e^{ik\lambda } \  \\
\times I_2 \Big (   \int_{j} ^{j+1} du\int_{j+|k|} ^{j+|k|+1} dv \vert u-v\vert ^{H-1}\int_{j} ^{j+1} du'(u'-y_{1})_+  ^{\frac{H}{2}-1} (u'-y_{2})_+  ^{\frac{H}{2}-1}\Big ).
\end{multline*}
But from the definition of a Rosenblatt process and the computation of $\int_{j} ^{j+1} du\int_{j+|k|} ^{j+|k|+1} dv\vert u-v\vert ^{H-1}$,
\begin{eqnarray*}
\widehat T'_{2, N}(g)=C^2 \,   \sqrt{ \frac{8(2H-1)}{H(1+H)^2}} \, \frac 1 { \pi \, N}\,  \sum_{|k|<N} \Big ( \int_{-\pi} ^{\pi } g(\lambda ) e^{ik\lambda } d\lambda \Big ) \big (|k+1|^{H+1}+ |k-1|^{H+1}-2|k|^{H+1}\big )  R_{N-|k|} ,
\end{eqnarray*}
where $(R_t)_t$ is a standard Rosenblatt process with parameter $H$.
But from Lemma \ref{l1}, with $\delta>0$ that can be chosen arbitrary small, we have $\int_{-\pi} ^{\pi } g(\lambda ) e^{ik\lambda } d\lambda =O((1+|k|)^{-2H+\delta})$  and since $\big (|k+1|^{H+1}+ |k-1|^{H+1}-2|k|^{H+1}\big )\sim \frac 1 2 (H+1)H \, k^{H-1}$ when $k$ is large enough, we deduce that
$$
 \sum_{k\in \Z} \Big ( \int_{-\pi} ^{\pi } g(\lambda ) e^{ik\lambda } d\lambda \Big ) \big (|k+1|^{H+1}+ |k-1|^{H+1}-2|k|^{H+1}\big ) <\infty.
$$
Moreover, since a Rosenblatt process is a $H$-self-similar process, $(N^{-H} R_{N-|k|})_{|k|<N} \simloi (R_{1-\frac {|k|}N})_{|k|<N}$. This is also a continuous process and therefore
\begin{eqnarray*}
N^{1-H} \widehat T'_{2, N}(g) &\simloi & 4C^2  \, \sqrt{ \frac{2(2H-1)}{H(1+H)^2}} \, \frac 1 { 2\pi}\,  \sum_{|k|<N} \Big ( \int_{-\pi} ^{\pi } g(\lambda ) e^{ik\lambda } d\lambda \Big ) \big (|k+1|^{H+1}+ |k-1|^{H+1}-2|k|^{H+1}\big )  R_{1-\frac {|k|} N} \\
&\limiteloiN &  4C^2  \, \sqrt{ \frac{2(2H-1)}{H(1+H)^2}} \,\int_{-\pi} ^{\pi } d\lambda g(\lambda ) \,  \frac 1 { 2\pi}\,  \sum_{k\in \Z } e^{ik\lambda }   \big (|k+1|^{H+1}+ |k-1|^{H+1}-2|k|^{H+1}\big )\,   R_{1}.
\end{eqnarray*}
But $\big (|k+1|^{H+1}+ |k-1|^{H+1}-2|k|^{H+1}\big )=2 \, r_{(H+1)/2,1}(k)$ (see (\ref{covY}))  and therefore from the definition of a spectral density,
\begin{eqnarray*}
N^{1-H} \widehat T_{2, N}(g)
&\limiteloiN &8  \,  \sqrt{ \frac{2(2H-1)}{H(1+H)^2}} \,\Big ( \int_{-\pi} ^{\pi }  g(\lambda )f_{(H+1)/2,C}(\lambda) \, d\lambda \Big )\,   R_{1}.
\end{eqnarray*}
\qed

\vskip0.2cm

  \begin{lemma}\label{A22}
Define the sequence of functions $(G_{N})_N$ by
\begin{multline*}
G_{N} (y_{1}, y_{2}) :=N ^{-H} \int_{-\pi } ^{\pi} g(\lambda) d\lambda \sum_{\vert k\vert <N} e^{-i\lambda k}  \sum_{j=0}^{N-|k|-1} \\
  \int_{j}^{j+1}du \int_{j+|k|} ^{j+|k|+1}dv \vert u-v\vert ^{H-1}  \Big [ (u-y_{1})_+^{\frac{H}{2}-1} (v-y_{2})_+^{\frac{H}{2}-1} - \int_j^{j+1}(u'-y_{1})_+^{\frac{H}{2}-1} (u'-y_{2})_+^{\frac{H}{2}-1}du' \Big ] .
\end{multline*}
Then,
$$
G_{N} (y_{1}, y_{2})  \limiteLL 0.
$$

\end{lemma}
{\bf Proof: }  1. On the one hand, define the sequence of functions $(G_{N,1})_N$ by
\begin{multline*}
G_{N,1} (y_{1}, y_{2}) :=N ^{-H} \int_{-\pi } ^{\pi} g(\lambda) d\lambda \sum_{\vert k\vert <N} e^{-i\lambda k}  \sum_{j=0}^{N-|k|-1} \\
 \int_{j}^{j+1}du \int_{j+|k|} ^{j+|k|+1}dv \vert u-v\vert ^{H-1}  (u-y_{1})_+^{\frac{H}{2}-1}\Big [ (v-y_{2})_+^{\frac{H}{2}-1} -  (u-y_{2})_+^{\frac{H}{2}-1} \Big ] .
\end{multline*}
Then we have
\begin{multline*}
\Vert G_{N,1} (y_{1}, y_{2}) \Vert ^{2} _{L^{2} (\mathbb{R}^{2})}= N ^{-2H}  \int_{\R^2}dy_1 \, dy_2 \, \sum_{\vert k_1\vert <N} \sum_{\vert k_2\vert <N}  \int_{-\pi } ^{\pi} g(\lambda_1)e^{-i\lambda_1 k_1} d\lambda_1
 \int_{-\pi } ^{\pi} g(\lambda_2)e^{-i\lambda_2 k_2} d\lambda_2  \\
 \times  \sum_{j_1=0}^{N-|k_1|-1}  \sum_{j_2=0}^{N-|k_2|-1} \Big ( \int_{j_1}^{j_1+1}du_1 \int_{j_1+|k_1|} ^{j_1+|k_1|+1}dv_1 \vert u_1-v_1\vert ^{H-1} (u_1-y_{1})^{\frac{H}{2}-1}_+\big [ (v_1-y_{2})_+^{\frac{H}{2}-1} -(u_1-y_{2})_+^{\frac{H}{2}-1} \big ] \Big ) \\
\times  \Big (  \int_{j_2}^{j_2+1}du_2 \int_{j_2+|k_2|} ^{j_2+|k_2|+1}dv_2 \vert u_2-v_2\vert ^{H-1} (u_2-y_{1})^{\frac{H}{2}-1}_+\big [ (v_2-y_{2})_+^{\frac{H}{2}-1} -(u_2-y_{2})_+^{\frac{H}{2}-1} \big ] \Big ).
\end{multline*}
Then,  using \eqref{i}, there exist $D_1>0$ and $D_2>0$ such as,
\begin{eqnarray}
\nonumber \Vert G_{N,1}\Vert ^{2} _{L^{2} (\mathbb{R}^{2})}&\hspace{-2mm}= &\hspace{-2mm} D_1 \,  N ^{-2H}  \sum_{\vert k_1\vert <N} \sum_{\vert k_2\vert <N}  \int_{-\pi } ^{\pi} g(\lambda_1)e^{-i\lambda_1 k_1} d\lambda_1
 \int_{-\pi } ^{\pi} g(\lambda_2)e^{-i\lambda_2 k_2} d\lambda_2   \\
\nonumber && \times   \sum_{j_1=0}^{N-|k_1|-1}  \sum_{j_2=0}^{N-|k_2|-1}\int_{j_1}^{j_1+1} \hspace{-3mm} du_1 \int_{j_1+|k_1|} ^{j_1+|k_1|+1}\hspace{-3mm}  dv_1 \int_{j_2}^{j_2+1} \hspace{-3mm}  du_2 \int_{j_2+|k_2|} ^{j_2+|k_2|+1}\hspace{-3mm}  dv_2 \, \vert u_1-v_1\vert ^{H-1} \vert u_2-v_2\vert ^{H-1} \vert u_1-u_2\vert ^{H-1} \\
\nonumber && \hspace{3cm} \times \big [ \vert u_1-u_2\vert ^{H-1}-\vert u_1-v_{2}\vert ^{H-1}-\vert u_2-v_{1}\vert ^{H-1}+ \vert v_{1}-v_{2}\vert ^{H-1}\big ]\\
\nonumber&&\hspace{-2.7cm}\leq \frac {D_2}  {N ^{2H}}  \sum_{\vert k_1\vert , \vert k_2\vert<N}  \frac 1 {[(1+|k_1|)(1+|k_2|) ]^{2H-\delta}} \hspace{-2mm}    \sum_{j_1,j_2=1} ^{N} \hspace{-1mm}    \int_{[0,1]^4}\hspace{-6mm} du_1  dv_1 du_2 dv_2 \,\big ( \vert u_1-v_1-|k_1|\vert\, \vert u_2-v_2-|k_2|\vert \, \vert u_1-u_2+j_1-j_2\vert \big ) ^{H-1} \\
\nonumber&& \hspace{-2.7cm}\times \Big | \vert u_1-u_2+j_1-j_2\vert ^{H-1}\hspace{-2mm}-\vert u_1-v_{2}+j_1-j_2-|k_2|\vert ^{H-1}\hspace{-2mm}-\vert u_2-v_{1}+j_2-j_1-|k_1|\vert ^{H-1}\hspace{-2mm}+ \vert v_{1}-v_{2}+j_1-j_2+|k_1|-|k_2|\vert ^{H-1}\Big | \\
\label{I} &  &\hspace{-2.7cm}\leq \frac {D_2}  {N ^{2H-1}}  \sum_{\vert k_1\vert <N} \sum_{\vert k_2\vert <N} \sum_{\vert j \vert <N} I(k_1,k_2,j)
\end{eqnarray}
with
\begin{multline}
 I(k_1,k_2,j) : =\frac 1 {[(1+|k_1|)(1+|k_2|) ]^{2H-\delta}}     \int_{[0,1]^4}\hspace{-6mm} du_1  dv_1 du_2 dv_2 \big ( \vert u_1-v_1-|k_1|\vert  \,\vert u_2-v_2-|k_2|\vert \, \vert u_1-u_2+j\vert \big )^{H-1} \\
\times \Big | \vert u_1-u_2+j\vert ^{H-1}-\vert u_1-v_{2}+j-|k_2|\vert ^{H-1}-\vert u_2-v_{1}-j-|k_1|\vert ^{H-1}+ \vert v_{1}-v_{2}+j+|k_1|-|k_2|\vert ^{H-1}\Big |~~
\end{multline}
and where both the last inequalities are obtained from changes of variables, Lemma \ref{l1} (with $\delta>0$ which can be chosen arbitrary small) and symmetry properties. For $|x|\geq 2$ and $(u,v)\in [0,1]^2$, $|x+u-v|^{H-1}\leq 2^{1-H} \, |x|^{H-1}$ and for $|x|< 2$ and $(u,v)\in [0,1]^2$, $u\neq v$, $|x+u-v|^{H-1}\leq |u-v|^{H-1}$. Moreover from a usual Taylor expansion of the function $x\mapsto |1+x|^{H-1}$, there exists $C>0$ such as
$$
\Big |\vert y+j\vert ^{H-1}-\vert j\vert ^{H-1}\big (1+(H-1)\frac y{|j|}\big )\Big |\leq C \, \frac y{|j|^{2-H}} \qquad \mbox{when $2|y|\leq j$}.
$$
Therefore, there exists $C>0$ such as
\begin{multline*}
\Big | \vert u_1-u_2+j\vert ^{H-1}-\vert u_1-v_{2}+j-|k_2|\vert ^{H-1}-\vert u_2-v_{1}-j-|k_1|\vert ^{H-1}+ \vert v_{1}-v_{2}+j+|k_1|-|k_2|\vert ^{H-1}\Big | \\
\leq C \, (1+|k_1|+|k_2|) \, |j|^{H-2} \\
\leq C \, \big [(1+|k_1|)^{2H-1+\delta}+(1+|k_2|)^{2H-1+\delta} \big ] \, |j|^{-H-\delta}
\end{multline*}
when $|j|> 2\max (|k_1|,|k_2|)$.
Then for $|j|> 2\max ( |k_1|,|k_2|)$ and for $|k_1|$ and $|k_2|$ large enough,
$$
 I(k_1,k_2,j) \leq \frac C {[(1+|k_1|)(1+|k_2|) ]^{1+H-\delta}}\,   \big [(1+|k_1|)^{2H-1+\delta}+(1+|k_2|)^{2H-1+\delta} \big ] \, \frac 1 {|j|^{1+\delta}}.
$$
Therefore,
\begin{equation} \label{MajI1}
\frac {D_2}  {N ^{2H-1}}  \sum_{\vert k_1\vert <N} \sum_{\vert k_2\vert <N} \sum_{|j|> 2\max (|k_1|,|k_2|)} I(k_1,k_2,j)  \limiteN 0.
\end{equation}
Now, for $2 \leq |j|\leq  2\max (|k_1|,|k_2|)$  and for $|k_1|$ and $|k_2|$ large enough,
\begin{multline*}
 I(k_1,k_2,j) \leq \frac C {[(1+|k_1|)(1+|k_2|) ]^{1+H-\delta}}     |j|^{H-1}  \\
\times  \Big ( \vert j\vert ^{H-1} + \vert u_1-v_{2}+j-|k_2|\vert ^{H-1}+\vert u_2-v_{1}-j-|k_1|\vert ^{H-1}+ \vert v_{1}-v_{2}+j+|k_1|-|k_2|\vert ^{H-1}\Big ).
\end{multline*}
But for $\beta\in \R^*$, since $1/2<H<1$, there exists $C>0$ such as for $|\beta| \leq M$,
$$
\sum_{j=2}^{M}  \frac 1 {(j(j+\beta))^{1-H}} \leq C\,  M^{2H-1}.
$$
Thus, for $2 \leq |j|\leq  2\max (|k_1|,|k_2|)$  and for $|k_1|$ and $|k_2|$ large enough,
$$
 I(k_1,k_2,j) \leq \frac C {[(1+|k_1|)(1+|k_2|) ]^{1+H-\delta}}\,\big ( \max (|k_1|,|k_2|) \big )^{2H-1}.
$$
Therefore,
\begin{equation} \label{MajI2}
\frac {D_2}  {N ^{2H-1}}  \sum_{\vert k_1\vert <N} \sum_{\vert k_2\vert <N} \sum_{2 \leq |j| \leq  2\max (|k_1|,|k_2|)} I(k_1,k_2,j)  \limiteN 0.
\end{equation}
We can easily add to this asymptotic behavior the cases $j=0$, $j=1$ or $j=-1$ thanks to the convergence of the integral defined in $[0,1]^4$. Finally, we obtain:
\begin{equation} \label{MajG1}
\Vert G_{N,1} (y_{1}, y_{2}) \Vert ^{2} _{L^{2} (\mathbb{R}^{2})} \limiteN 0.
\end{equation}
2. On the other hand, define the sequence of functions $(G_{N,2})_N$ by
\begin{multline*}
G_{N,2} (y_{1}, y_{2}) :=N ^{-H} \int_{-\pi } ^{\pi} g(\lambda) d\lambda \sum_{\vert k\vert <N} e^{-i\lambda k}  \sum_{j=0}^{N-|k|-1} \\
 \int_{j}^{j+1}du \int_{j+|k|} ^{j+|k|+1}dv \int_{j}^{j+1}du'  \vert u-v\vert ^{H-1} \Big [  (u-y_{1})_+^{\frac{H}{2}-1}(u-y_{2})_+^{\frac{H}{2}-1} -(u'-y_{1})_+^{\frac{H}{2}-1}(u'-y_{2})_+^{\frac{H}{2}-1}\Big ] .
\end{multline*}
Then following the same kind of computations and expansions than for $\Vert G_{N,1} (y_{1}, y_{2}) \Vert ^{2} _{L^{2} (\mathbb{R}^{2})}$, we have
\begin{eqnarray*}
\Vert G_{N,2} (y_{1}, y_{2}) \Vert ^{2} _{L^{2} (\mathbb{R}^{2})}&\hspace{-3mm} =& \hspace{-3mm} N ^{-2H}  \sum_{\vert k_1\vert <N} \sum_{\vert k_2\vert <N}   \sum_{j_1=0}^{N-|k_1|-1}   \sum_{j_2=0}^{N-|k_2|-1}  \int_{-\pi } ^{\pi} g(\lambda_1)e^{-i\lambda_1 k_1} d\lambda_1
 \int_{-\pi } ^{\pi} g(\lambda_2)e^{-i\lambda_2 k_2} d\lambda_2  \\
&\hspace{-3mm}\times & \hspace{-3mm}\hspace{-0cm}   \int_{j_1}^{j_1+1} \hspace{-5mm}du_1 \int_{j_1+|k_1|} ^{j_1+|k_1|+1}\hspace{-7mm}dv_1 \int_{j_1}^{j_1+1} \hspace{-5mm}du'_1  \int_{j_2}^{j_2+1}\hspace{-5mm}du_2 \int_{j_2+|k_2|} ^{j_2+|k_2|+1}\hspace{-7mm}dv_2  \int_{j_2}^{j_2+1}\hspace{-5mm}du'_2  \, |u_1-v_1|^{H-1} |u_2-v_2|^{H-1} \\
&& \hspace{+2.5cm}  \times
\Big ( |u_1-u_{2}|^{2H-2} +|u'_1-u'_{2}|^{2H-2}-|u'_2-u_{1}|^{2H-2} -|u'_1-u'_{2}|^{2H-2}\Big )\\
&\hspace{-3mm} \leq &\hspace{-3mm} \frac {C}  {N ^{2H-1}} \hspace{-0.4cm}  \sum_{\vert k_1\vert, \vert k_2\vert <N} \frac 1 {[(1+|k_1|)(1+|k_2|) ]^{2H-\delta}}\hspace{-1mm} \int_{[0,1]^6}\hspace{-6mm}  du_1 dv_1 du'_1 du_2 dv_2 du'_2 \big (|u_1-v_1+|k_1||\,|u_2-v_2+|k_2||\big )^{H-1}  \\
&& \times \sum_{|j|<N}   \Big | |u_1-u_{2}+j|^{2H-2}  \hspace{-1mm}-|u'_1+j-u_{2}|^{2H-2} \hspace{-1mm}-|u_{1}+j-u'_2|^{2H-2} \hspace{-1mm}+|j+u'_1-u'_2|^{2H-2}\Big |   \\
&\hspace{-3mm}\leq &\hspace{-3mm} \frac {C}  {N ^{2H-1}}  \sum_{\vert k_1\vert <N} \sum_{\vert k_2\vert <N} \sum_{|j|<N}  \hspace{-1mm}   \frac 1 {[(1+|k_1|)(1+|k_2|) ]^{H+1-\delta}} \, (1+|j|)^{2H-3}.
\end{eqnarray*}
As a consequence, since all the previous sums are finite, we deduce that
\begin{equation} \label{MajG2}
\Vert G_{N,2} (y_{1}, y_{2}) \Vert ^{2} _{L^{2} (\mathbb{R}^{2})} \limiteN 0.
\end{equation}
Then \eqref{MajG1} and \eqref{MajG2} imply Lemma \ref{A22} holds.  \qed
~\\

{\bf Proof of Proposition \ref{p3}: } In a first step we use the isometry of multiple integrals (\ref{isom}) and the fact that the $L^{2}$ norm of the symmetrized function is less than the $L^{2}$ norm of the function itself. In a second step we change the order of integration and we use (\ref{i}). Then, for the last bound below, we use  the same lines as in  the proof of Lemma \ref{A22}. We have
\begin{eqnarray*}
\E \left( N^{1-H} \widehat T_{4, N}(g) \right) ^{2} &\leq  & N^{-2H} \int_{\mathbb{R}^4 }  dy_{1} dy_{2}  dy_{3} dy_{4}  \int_{[-\pi,\pi]^2}   d\lambda _{1} d\lambda _{2} g(\lambda _{1}) g(\lambda _{2} ) \sum_{|k_{1}|<N}\sum_{|k_{2}|<N} e^{-ik_{1}\lambda _{1}}e^{-ik_{2}\lambda _{2}}  \\
&& \sum_{j_{1}=0}^{N-|k_1|-1} \sum_{j_{2}=0}^{N-|k_2|-1}  \int_{j_{1}} ^{j_{1}+1} du \int_{j_{1}+|k_{1}|} ^{j_{1}+|k_{1}|+1}dv (u-y_{1})_{+} ^{\frac{H}{2}-1} (u-y_{2})_{+} ^{\frac{H}{2}-1} (v-y_{3})_{+} ^{\frac{H}{2}-1}(v-y_{4})_{+} ^{\frac{H}{2}-1}\\
&& \times \int_{j_{2}} ^{j_{2}+1} du' \int_{j_{2}+|k_{2}|} ^{j_{2}+|k_{2}|+1}dv' (u'-y_{1})_{+} ^{\frac{H}{2}-1} (u'-y_{2})_{+} ^{\frac{H}{2}-1} (v-y_{3})_{+} ^{\frac{H}{2}-1}(v'-y_{4})_{+} ^{\frac{H}{2}-1} \\
&=&\beta^4(1-H,\frac H 2) \, N^{-2H}\sum_{|k_{1}|<N}\sum_{|k_{2}|<N} \big (\int_{-\pi}^\pi  e^{-ik_{1}\lambda _{1}} g(\lambda _{1})  d\lambda _{1}\big ) \big (\int_{-\pi}^\pi  e^{-ik_{2}\lambda _{2}} g(\lambda _{2})  d\lambda _{2 }\big )\\
&&\times  \sum_{j_{1}=0}^{N-|k_1|-1} \sum_{j_{2}=0}^{N-|k_2|-1}  \int_{j_{1}} ^{j_{1}+1}  \int_{j_{1}+|k_{1}|} ^{j_{1}+|k_{1}|+1} \int_{j_{2}} ^{j_{2}+1}  \int_{j_{2}+|k_{2}|} ^{j_{2}+|k_{2}|+1}\hspace{-5mm} du\, dv\,  du' \,  dv'\,  \vert u-u'\vert ^{2H-2}\vert v-v'\vert ^{2H-2} \\
&\leq  &  c\, N^{1-2H}\sum_{|k_{1}|<N}\sum_{|k_{2}|<N}  \frac 1 {\big ((1+|k_{1}|)(1+|k_{2}|)\big ) ^{2H-\delta }} \sum_{|j|<N}\big | r_{H,1}(j) \, r_{H,1}(j+|k_1|-|k_2|) \big |.
\end{eqnarray*}
As in the proof of Lemma \ref{A22}, we decompose the previous right left hand term in two parts: firstly, when $|j|>2\max(|k_1|,|k_2|)$, we have $r_{H,1}(j) \leq c \, (1+|j|)^{2H-2}$ and $r_{H,1}(j+|k_1|-|k_2|) \leq c \, (1+|j|)^{2H-2}$. Then, with $c>0$, 
\begin{eqnarray}
\nonumber  && N^{1-2H} \sum_{|k_{1}|<N}\sum_{|k_{2}|<N}  \frac 1 {\big ((1+|k_{1}|)(1+|k_{2}|)\big ) ^{2H-\delta }} \sum_{2\max(|k_1|,|k_2|)<|j|<N}\big |  r_{H,1}(j) \, r_{H,1}(j+|k_1|-|k_2|) \big |  \\
\nonumber & & \hspace{1cm}\leq  c \,  N^{1-2H} \sum_{|k_{1}|<N}\sum_{|k_{2}|<N}  \frac 1 {\big ((1+|k_{1}|)(1+|k_{2}|)\big ) ^{2H-\delta }}  \, N^{4H-3} \\
\label{casj1}& & \hspace{1cm}\leq  c \,  N^{2H-2 } ,
\end{eqnarray}
since $2H-\delta>0$ because $H>1/2$ and $\delta>0$ can be chosen arbitrary small. \\
Secondly, when $|j|\leq \max(|k_1|,|k_2|)$, using $\big | r_{H,1}(j+|k_1|-|k_2|)\big |  \leq 1$, we obtain
\begin{eqnarray}
\nonumber  && N^{1-2H} \sum_{|k_{1}|<N}\sum_{|k_{2}|<N}  \frac 1 {\big ((1+|k_{1}|)(1+|k_{2}|)\big ) ^{2H-\delta }} \sum_{|j|\leq 2\max(|k_1|,|k_2|)} \big | r_{H,1}(j) \, r_{H,1}(j+|k_1|-|k_2|)\big |  \\
\nonumber & & \hspace{1cm}\leq  c \,  N^{1-2H} \sum_{|k_{1}|<N}\sum_{|k_{2}|<N}  \frac 1 {\big ((1+|k_{1}|)(1+|k_{2}|)\big ) ^{2H-\delta }}  \, \max(|k_1|,|k_2|)^{2H-1} \\
\label{casj2}& & \hspace{1cm}\leq   c \,  N^{1-2H+\delta}.
\end{eqnarray}
As a consequence, from \eqref{casj1} and \eqref{casj2},
$$
\E \left( N^{1-H} \widehat T_{4, N}(g) \right) ^{2} \limiteN 0,
$$
and the conclusion of Proposition \ref{p3} follows.
\qed
~\\
~\\
{\bf Proof of Proposition \ref{A4}: }The proof follows the lines of end of the proof of Theorem 2 in \cite{FT}, p. 528-529. Notice that
$$\E \widehat J_{N}\Big (\frac{\partial}{ \partial H}\big (\frac{1}{g_{H}}\big )\Big )= \sum _{\vert k\vert \leq N} \int_{-\pi } ^{\pi } d\lambda \frac{\partial}{ \partial H}\Big ( \frac{1}{g_{H} (\lambda)} \Big )\frac{N-k}{N} r_{H,C}(k) e^{-ik\lambda  } .
$$
Denote
$$
e_{k}= \int_{-\pi } ^{\pi } e^{ik\lambda }   \frac{\partial}{ \partial H}\Big ( \frac{1}{g_{H} (\lambda)}\Big ) \, d\lambda.
$$
Then
$$
\E \widehat J_{N}\Big (\frac{\partial}{ \partial H}\big (\frac{1}{g_{H}}\big )\Big )= \frac {1}{(2\pi)^2 N} \,  \sum_{j=1}^N \sum_{k=1}^N r_{H,C}(j-k) \, e_{j-k}
$$
From (3.3) in \cite{FT}, we know that
\begin{equation}\label{h(0)}
\sum_{k\in \Z} r_{H,C}(k)\, e_k=0.
\end{equation}
Therefore,
\begin{eqnarray}
\nonumber N^{1-H} \,  \E \widehat J_{N}\Big (\frac{\partial}{ \partial H}\big (\frac{1}{g_{H}}\big )\Big )&=& N^{1-H} \sum_{|k|<N} \big ( 1 - \frac k N \big ) \, r_{H,C}(k)\, e_k \\
\label{somme}&=& N^{1-H} \sum_{|k|<N}  r_{H,C}(k)\, e_k -N^{-H} \sum_{|k|<N} k\, r_{H,C}(k)\, e_k.
\end{eqnarray}
From Lemma \ref{l1} and since $r_{H,C}(k)=O(|k|^{2H-2}),~ |k|\to \infty$, we have for all $\delta>0$
\begin{equation}\label{er}
 r_{H,C}(k)\, e_k=O\big (|k|^{-2+\delta}\big ).
\end{equation}
Then, from a usual comparison with a Riemann integral, $\sum_{|k|<N} k\, r_{H,C}(k)\, e_k=O(N^\delta)$. Moreover, from \eqref{h(0)}, $\sum_{|k|<N}  r_{H,C}(k)\, e_k=-\sum_{|k|\geq N}  r_{H,C}(k)\, e_k$ and also from a usual comparison with a Riemann integral and \eqref{er}, we know that $\sum_{|k| \geq N}  r_{H,C}(k)\, e_k = O(N^{-1+\delta})$. As a consequence, from \eqref{somme},
 \begin{eqnarray*}
\nonumber N^{1-H} \,  \E \widehat J_{N}\Big (\frac{\partial}{ \partial H}\big (\frac{1}{g_{H}}\big )\Big )&=& O\big (N^{\delta-H} \big).
\end{eqnarray*}
As $\delta$  can be chosen arbitrary small, we deduce that Proposition \ref{A4} holds.
\qed
~\\
~\\
{\bf Proof of Proposition \ref{p4}: } Note that the Fourier coefficients of $\widehat I_{N}$ are given by
\begin{equation*}
\int_{-\pi} ^{\pi} e^{ikx } \widehat  I_{N}(x) dx = \widehat{r}_{N} (k) 1_{(\vert k\vert <N)}.
\end{equation*} The proofs of Proposition \ref{p2} and \ref{p3} imply that
\begin{equation}\label{as}
 \widehat{r} _{N} (k) \limiteL2 r_{H,C}(k).
\end{equation}
 The above convergence holds almost surely by an argument in \cite{D} (Theorem 7.1, p. 493).  \qed
~\\
~\\
{\bf Proof of Theorem \ref{asymptoH}: } For establishing the strong convergence, define for $h \in (1/2,1)$,
\begin{equation} \label{VN}
\widehat V_N(h)=\frac 1 {2\pi} \, \widehat J_N(1/g_h). 
\end{equation}
From Proposition \ref{p4}, $\widehat V_N(h) \limiteasN V(h)=\sigma^2 \int_{-\pi}^\pi g_H(\lambda)/g_h(\lambda)d\lambda$. As $H=\mbox{Argmin}_{h\in (1/2,1)} V(h)$ and  $\widehat H_N=\mbox{Argmin}_{h\in (1/2,1)}\widehat V_N(h)$, from usual arguments (see for instance \cite{FT}), then $\widehat H_N \limiteasN H$. \\
~\\
For proving the non-central theorem, apply Proposition \ref{p7} to $g=\frac{ \partial  }{\partial H } \left( \frac{1}{g_{H} } \right)$ and from Proposition \ref{A4}, we obtain:
$$
N^{1-H}\,  \frac{\partial}{\partial h}\widehat V_N(H) \limiteloiN \frac 1 {2\pi} \, \sqrt{ \frac{8(2H-1)}{H(1-H)^2}} \,\Big ( \int_{-\pi} ^{\pi } \frac{\partial}{\partial H}  \Big (\frac {1}{g_H(\lambda)}\Big ) \,g_{(H+1)/2}(\lambda) \,  d\lambda \Big ) \, R,
$$
with $R$ a standard Rosenblatt random variable. Now use the following result established in Lemma 2 in \cite{FT}.
Suppose that $(b_{N})_N$ is a sequence of real numbers such that $b_{N} \limiteN +\infty$. Assume that
\begin{equation}\label{bN}
b_{N} \frac{\partial}{\partial h} \widehat V_N(H) \limiteloiN Y,
\end{equation}
 where $Y$ is a random variable. Then
\begin{equation*}
b_{N} (\widehat{H} _{N} -H) \limiteloiN  - 2\pi\, \Big (\int_{-\pi} ^{\pi} f_{H,C} (\lambda) \frac{ \partial ^{2} }{\partial H^{2} } \left( \frac{1}{g_{H} (\lambda)} \right) d\lambda \Big ) ^{-1} \, Y.
\end{equation*}
Using $b_N=N^{1-H}$, $Y=8  \,  \sqrt{ \frac{2(2H-1)}{H(1-H)^2}} \,\Big ( \int_{-\pi} ^{\pi }  g(\lambda )f_{(H+1)/2,C}(\lambda) \, d\lambda \Big )\,   R_{1} $,  Proposition \ref{p7} applied to $ \frac{\partial}{\partial H} \big (\frac{1}{g_{H}} \big )$ and Proposition \ref{A4}, then \eqref{bN} holds. This implies Theorem \ref{asymptoH}. \qed
~\\
~\\
{\bf Proof of Corollary \ref{sigma}: } Using $\widehat V_N(h)$ defined in \eqref{VN}, it is sufficient to write 
$$
\widehat C_N=\big (\mu(\widehat H_N) \, \widehat V_N(\widehat H_N) \big )^{1/2}, 
$$
with 
\begin{equation} \label{eqbeta}
\mu(h)=2\pi \, \frac {a_{h}  }{h \, \Gamma(2h)\sin(\pi h)} \qquad \mbox{for $h\in (1/2,1)$}, 
\end{equation}
and the strong consistencies established in Proposition \ref{p4} and Theorem \ref{asymptoH} allow to show the almost sure convergence of $\widehat C_N$. \\
~\\
For establishing the non-central limit theorem, the Taylor's formula implies that
$$
\widehat V_N(H)=\widehat V_N(\widehat H_N)+\frac 1 2 \, (H-\widehat H_N)^2
 \Big (\frac {\partial^2 } {\partial
h^2}\widehat V_N(\underline{H}_N) \Big ),
$$
with probability $1$, and with $| \underline{H}_N- H|< |\widehat{H}_N - H |$. From previous Theorem
\ref{asymptoH}, it follows
\begin{equation}\label{eq1}
N^{1-H}( \widehat V_N(H)- \sigma^2)
=N^{1-H}( \widehat V_N(\widehat{H}_N)-\sigma^2)+ {\cal O}_p(N^{H-1}).
\end{equation}
Moreover, from Proposition \ref{p7}, 
\begin{equation}\label{eq2}
N^{1-H}\big ( \widehat V_N(H)- \E \big ( \widehat V_N(H)\big )\big ) \limiteloiN \frac 4 {\pi} \,C^2 \,  \beta(H) \,  R.
\end{equation}
with $\beta(H)=\sqrt{ \frac{2(2H-1)}{H(1-H)^2}} \,\Big ( \int_{-\pi} ^{\pi }  g^{-1}_H(\lambda )f_{(H+1)/2,1}(\lambda) \, d\lambda \Big ) $. 
Since $ \E \big ( \widehat V_N(H)\big )=\sigma^2 + o(N^{-1/2})$ (see for instance \cite{FT}), from \eqref{eq1} and \eqref{eq2}, we deduce: 
\begin{equation}\label{eq3}
N^{1-H}( \widehat V_N(\widehat{H}_N)-\sigma^2) \limiteloiN
\frac 4 {\pi} \, \beta(H) \, R.
\end{equation}
From Theorem \ref{asymptoH} and using the Delta-Method, 
\begin{equation}\label{eq4}
 N^{1-H} \big ( \mu(\widehat{H}_N) - \mu(H)\big )  \limiteloiN  (\mu'(H))^2  \, \gamma(H) \, R.
\end{equation}
Thus, using \eqref{eq3} and \eqref{eq4}, we obtain:
\begin{equation}\label{eq5}
 N^{1-H} \big ( \mu(\widehat{H}_N) \widehat V_N(\widehat{H}_N) - \mu(H)\sigma^2 )  \limiteloiN  C^2 \, \Big (  (\mu'(H))^2 \, \gamma(H) \,\frac 1 {\mu(H)} +\frac 4 {\pi} \,\beta(H) \,  \mu(H) \Big ) \,   R.
\end{equation}
As a consequence, using again the Delta-Method, we have:
\begin{equation}\label{eq6}
 N^{1-H} \big (\widehat C_N- C)  \limiteloiN  \frac C {4 \, \mu(H)}\,  \Big (  (\mu'(H))^2 \, \gamma(H)  +\frac 4 {\pi} \,\beta(H) \,  \mu^2(H) \Big ) \,   R,
\end{equation}
achieving the proof of Corollary \ref{sigma}.

\qed

\end{document}